\documentclass{amsart}
\usepackage{amsthm}
\usepackage{amstext}
\usepackage{amssymb}
\usepackage{hyperref}

\makeatletter
\numberwithin{equation}{section}
\numberwithin{figure}{section}
\theoremstyle{plain}
\newtheorem{thm}{Theorem}[section]

\allowdisplaybreaks
\usepackage[sort&compress,square,numbers]{natbib}

\newcommand\relphantom[1]{\mathrel{\phantom{#1}}}
\makeatother

\begin{document}

\title[Symmetry for $q$-Euler polynomials]{Some identities of symmetry for $q$-Euler polynomials under the
symmetric group of degree $n$ arising from fermionic $p$-adic $q$-integrals
on $\mathbb{Z}_{p}$}
\author{Dmitry V. Dolgy}
\address{Institute of Natural Sciences, Far Eastern Federal University, 690950 Vladivostok
 Russia}
\email{$d_{-}dol@mail.ru$}

\author{Dae San Kim}
\address{Department of Mathematics, Sogang University, Seoul 121-742, Republic
of Korea}
\email{dskim@sogang.ac.kr}

\author{Taekyun Kim}
\address{Department of Mathematics, Kwangwoon University, Seoul 139-701, Republic
of Korea}
\email{tkkim@kw.ac.kr}


\begin{abstract}
In this paper, we investigate some new symmetric identities for the
$q$-Euler polynomials under the symmetric group of degree $n$ which
are derived from fermionic $p$-adic $q$-integrals on $\mathbb{Z}_{p}$.
\end{abstract}

\maketitle
\global\long\def\acl#1#2{\left\langle \left.#1\right|#2\right\rangle }

\global\long\def\acr#1#2{\left\langle #1\left|#2\right.\right\rangle }

\global\long\def\Li{\mathrm{Li}}

\global\long\def\Zp{\mathbb{Z}_{p}}

\section{Introduction}

Let $p$ be a fixed prime number such that $p\equiv1\pmod{2}$. Throughout
this paper, $\mathbb{Z}_{p}$, $\mathbb{Q}_{p}$ and $\mathbb{C}_{p}$
will denote the ring of $p$-adic integers, the field of $p$-adic
rational numbers and the completion of the algebraic closure of $\mathbb{Q}_{p}$.
Let $q$ be an indeterminate in $\mathbb{C}_{p}$ such that $\left|1-q\right|_{p}<p^{-\frac{1}{p-1}}$.
The $p$-adic norm is normalized as $\left|p\right|_{p}=\frac{1}{p}$
and the $q$-analogue of the number $x$ is defined as $\left[x\right]_{q}=\frac{1-q^{x}}{1-q}$.
Note that $\lim_{q\rightarrow1}\left[x\right]_{q}=x$.

As is well known, the Euler numbers are defined by
\[
E_{0}=1,\quad\left(E+1\right)^{n}+E_{n}=2\delta_{0,n},\quad\left(n\in\mathbb{N}\cup\left\{ 0\right\} \right),
\]
with the usual convention about replacing $E^{n}$ by $E_{n}$(see
\cite{key-1,key-2,key-3,key-4,key-5,key-6,key-7,key-8,key-9,key-10,key-11,key-12,key-13,key-14}).

The Euler polynomials are given by
\[
E_{n}\left(x\right)=\sum_{l=0}^{n}\binom{n}{l}x^{n-l}E_{l}=\left(E+x\right)^{n},\quad\left(n\ge0\right),\quad\left(\text{see \cite{key-1,key-2}}\right).
\]
In \cite{key-5}, Kim introduced Carlitz-type $q$-Euler numbers as follows:
\begin{equation}
\mathcal{E}_{0,q}=1,\quad q\left(q\mathcal{E}_{q}+1\right)^{n}+\mathcal{E}_{n,q}=\left[2\right]_{q}\delta_{0,n},\quad\left(n\ge0\right),\quad\left(\text{see \cite{key-5}}\right),\label{eq:1}
\end{equation}
with the usual convention about replacing $\mathcal{E}_{q}^{n}$ by
$\mathcal{E}_{n,q}$.

The Carlitz-type $q$-Euler polynomials are also defined as
\begin{equation}
\mathcal{E}_{n,q}\left(x\right)=\left(q^{x}\mathcal{E}_{q}+\left[x\right]_{q}\right)^{n}=\sum_{l=0}^{n}\binom{n}{l}q^{lx}\mathcal{E}_{l,q}\left[x\right]_{q}^{n-l},\quad\left(\text{see \cite{key-3,key-5}}\right).\label{eq:2}
\end{equation}

Let $C\left(\Zp\right)$ be the space of all $\mathbb{C}_{p}$-valued continuous functions
on $\Zp$. Then, for $f\in C\left(\Zp\right)$, the fermionic $p$-adic
$q$-integral on $\Zp$ is defined by Kim as
\begin{align}
I_{-q}\left(f\right) & =\int_{\Zp}f\left(x\right)d\mu_{-q}\left(x\right)\label{eq:3}\\
 & =\lim_{N\rightarrow\infty}\frac{1}{\left[p^{N}\right]_{-q}}\sum_{x=0}^{p^{N}-1}f\left(x\right)\left(-q\right)^{x}\nonumber \\
 & =\lim_{N\rightarrow\infty}\frac{1+q}{1+q^{p^{N}}}\sum_{x=0}^{p^{N}-1}f\left(x\right)\left(-q\right)^{x},\quad\left(\text{see \cite{key-5,key-6,key-7,key-8,key-9,key-10,key-11}}\right).\nonumber
\end{align}

From (\ref{eq:3}), we note that
\begin{equation}
q^{n}I_{-q}\left(f_{n}\right)+\left(-1\right)^{n-1}I_{-q}\left(f\right)=\left[2\right]_{q}\sum_{l=0}^{n-1}\left(-1\right)^{n-1-l}f\left(l\right),\quad\left(n\in\mathbb{N}\right),\quad\left(\text{see \cite{key-5}}\right).\label{eq:4}
\end{equation}

The Carlitz-type $q$-Euler polynomials can be represented by the
fermionic $p$-adic $q$-integral on $\Zp$ as follows:
\begin{equation}
\mathcal{E}_{n,q}\left(x\right)=\int_{\Zp}\left[x+y\right]_{q}^{n}d\mu_{-q}\left(y\right),\quad\left(n\ge0\right),\quad\left(\text{see \cite{key-5}}\right).\label{eq:5}
\end{equation}

Thus, by (\ref{eq:5}), we get
\begin{align}
 & \mathcal{E}_{n,q}\left(x\right)\label{eq:6}\\
= & \sum_{l=0}^{n}\binom{n}{l}q^{lx}\int_{\Zp}\left[y\right]_{q}^{l}d\mu_{q}\left(y\right)\left[x\right]_{q}^{n-l}\nonumber \\
= & \sum_{l=0}^{n}\binom{n}{l}q^{lx}\mathcal{E}_{l,q}\left[x\right]_{q}^{n-l},\quad\left(\text{see \cite{key-5}}\right).\nonumber
\end{align}

From (\ref{eq:4}), we can easily derive
\begin{equation}
q\int_{\Zp}\left[x+1\right]_{q}^{n}d\mu_{-q}\left(x\right)+\int_{\Zp}\left[x\right]_{q}^{n}d\mu_{-q}\left(x\right)=\left[2\right]_{q}\delta_{0,n},\quad\left(n\in\mathbb{N}\cup\left\{ 0\right\} \right).\label{eq:7}
\end{equation}

The equation (\ref{eq:7}) is equivalent to
\begin{equation}
q\mathcal{E}_{n.q}\left(1\right)+\mathcal{E}_{n,q}=\left[2\right]_{q}\delta_{0,n},\quad\left(n\geq0\right).
\label{eq:8}
\end{equation}

The purpose of this paper is to give some new symmetric identities
for the Carlitz-type $q$-Euler polynomials under the symmetric
group of degree $n$ which are derived from fermionic $p$-adic $q$-integrals
on $\Zp$.

\section{Symmetric identities for $\mathcal{E}_{n,q}\left(x\right)$ under
$S_{n}$}

Let $w_{1},w_{2},\dots,w_{n}\in\mathbb{N}$ such that $w_{1}\equiv w_{2}\equiv w_{3}\equiv\cdots\equiv w_{n}\equiv1\pmod{2}$.
Then, we have
\begin{align}
 &\relphantom{=} \int_{\Zp}e^{\left[\left(\prod_{j=1}^{n-1}w_{j}\right)y+\left(\prod_{j=1}^{n}w_{j}\right)x+w_{n}\sum_{j=1}^{n-1}\left(\prod_{\substack{i=1\\
i\neq j
}
}^{n-1}w_{i}\right)k_{j}\right]_{q}t}d\mu_{-q^{w_{1}\cdots w_{n-1}}}\left(y\right)\label{eq:8}\\
 &= \lim_{N\rightarrow\infty}\frac{1}{\left[p^{N}\right]_{q^{w_{1}\cdots w_{n-1}}}}\nonumber \\
&\relphantom{=}\times \sum_{y=0}^{p^{N}-1}e^{\left[\left(\prod_{j=1}^{n-1}w_{j}\right)y+\left(\prod_{j=1}^{n}w_{j}\right)x+w_{n}\sum_{j=1}^{n-1}\left(\prod_{\substack{i=1\\
i\neq j
}
}^{n-1}w_{i}\right)k_{j}\right]_{q}t}\left(-q^{w_{1}\cdots w_{n-1}}\right)^{y}\nonumber\\
= & \frac{1}{2}\lim_{N\rightarrow\infty}\left[2\right]_{q^{w_{1}\cdots w_{n-1}}}\nonumber\\
&\times\sum_{m=0}^{w_{n}-1}\sum_{y=0}^{p^{N}-1}e^{\left[\left(\prod_{j=1}^{n-1}w_{j}\right)\left(m+w_{n}y\right)+\left(\prod_{j=1}^{n}w_{j}\right)x+w_{n}\sum_{j=1}^{n}\left(\prod_{\substack{i=1\\
i\neq j
}
}^{n}w_{i}\right)k_{j}\right]_{q}t}\nonumber \\
 & \times\left(-1\right)^{m+y}q^{w_{1}\cdots w_{n-1}\left(m+w_{n}y\right)}.\nonumber
\end{align}

Thus, by (\ref{eq:8}), we get
\begin{align}
 &\relphantom{=} \frac{1}{\left[2\right]_{q^{w_{1}\cdots w_{n-1}}}}\prod_{l=1}^{n-1}\sum_{k_{l}=0}^{w_{l}-1}\left(-1\right)^{\sum_{i=1}^{n-1}k_{i}}q^{w_{n}\sum_{j=1}^{n-1}\left(\prod_{\substack{i=1\\
i\neq j
}
}^{n-1}w_{i}\right)k_{j}}\label{eq:9}\\
&\relphantom{=}\times\int_{\Zp}e^{\left[\left(\prod_{j=1}^{n-1}w_{j}\right)y+\left(\prod_{j=1}^{n}w_{j}\right)x+w_{n}\sum_{j=1}^{n-1}\left(\prod_{\substack{i=1\\
i\neq j
}
}^{n-1}w_{i}\right)k_{j}\right]_{q}t}d\mu_{-q^{w_{1}\cdots w_{n-1}}}\left(y\right)\nonumber\\
&=  \frac{1}{2}\lim_{N\rightarrow\infty}\prod_{l=1}^{n-1}\sum_{k_{l}=0}^{w_{l}-1}\sum_{m=0}^{w_{n}-1}\sum_{y=0}^{p^{N}-1}\left(-1\right)^{\sum_{i=1}^{n-1}k_{i}+m+y}\nonumber\\
&\relphantom{=}\times q^{w_{n}\sum_{j=1}^{n-1}\left(\prod_{\substack{i=1\\
i\neq j
}
}^{n-1}w_{i}\right)k_{j}+\left(\prod_{j=1}^{n-1}w_{j}\right)m+\left(\prod_{j=1}^{n}w_{j}\right)y}\nonumber \\
 &\relphantom{=} \times e^{\left[\left(\prod_{j=1}^{n-1}w_{j}\right)\left(m+w_{n}y\right)+\left(\prod_{j=1}^{n}w_{j}\right)x+w_{n}\sum_{j=1}^{n-1}\left(\prod_{\substack{i=1\\
i\neq j
}
}^{n-1}w_{i}\right)k_{j}\right]_{q}t}.\nonumber
\end{align}

As this expression is invariant under any permutation $\sigma\in S_{n}$,
we have the following theorem.
\begin{thm}
\label{thm:1}Let $w_{1},w_{2},\dots,w_{n}\in\mathbb{N}$ such that
$w_{1}\equiv w_{2}\equiv\cdots\equiv w_{n}\equiv1\pmod{2}$. Then,
the following expressions
\begin{align*}
 & \frac{1}{\left[2\right]_{q^{w_{\sigma\left(1\right)}\cdots w_{\sigma\left(n-1\right)}}}}\prod_{l=1}^{n-1}\sum_{k_{l}=0}^{w_{\sigma\left(l\right)}-1}\left(-1\right)^{\sum_{i=1}^{n-1}k_{i}}q^{w_{\sigma\left(n\right)}\sum_{j=1}^{n-1}\left(\prod_{\substack{i=1\\
i\neq j
}
}^{n-1}w_{\sigma\left(i\right)}\right)k_{j}}\\
&\times\int_{\Zp}e^{\left[\left(\prod_{j=1}^{n-1}w_{\sigma\left(j\right)}\right)y+\left(\prod_{j=1}^{n}w_{j}\right)x+w_{\sigma\left(n\right)}\sum_{j=1}^{n-1}\left(\prod_{\substack{i=1\\
i\neq j
}
}^{n-1}w_{\sigma\left(i\right)}\right)k_{j}\right]_{q}t}d\mu_{q^{w_{\sigma\left(1\right)}\cdots w_{\sigma\left(n-1\right)}}}\left(y\right)
\end{align*}
are the same for any $\sigma\in S_{n}$, $\left(n\ge1\right).$
\end{thm}
Now, we observe that
\begin{align}
 & \left[\left(\prod_{j=1}^{n-1}w_{j}\right)y+\left(\prod_{j=1}^{n}w_{j}\right)x+w_{n}\sum_{j=1}^{n-1}\left(\prod_{\substack{i=1\\
i\neq j
}
}^{n-1}w_{j}\right)k_{j}\right]_{q}t\label{eq:10}\\
= & \left[\prod_{j=1}^{n-1}w_{j}\right]_{q}\left[y+w_{n}x+w_{n}\sum_{j=1}^{n-1}\frac{k_{j}}{w_{j}}\right]_{q^{w_{1}\cdots w_{n-1}}}.\nonumber
\end{align}

By (\ref{eq:10}), we get
\begin{align}
 & \int_{\Zp}e^{\left[\left(\prod_{j=1}^{n-1}w_{j}\right)y+\left(\prod_{j=1}^{n}w_{j}\right)x+w_{n}\sum_{j=1}^{n-1}\left(\prod_{\substack{i=1\\
i\neq j
}
}^{n-1}w_{i}\right)k_{j}\right]_{q}t}d\mu_{-q^{w_{1}\cdots w_{n-1}}}\left(y\right)\label{eq:11}\\
= & \sum_{m=0}^{\infty}\left[\prod_{j=1}^{n-1}w_{j}\right]_{q}^{m}\int_{\Zp}\left[y+w_{n}x+w_{n}\sum_{j=1}^{n-1}\frac{k_{j}}{w_{j}}\right]_{q^{w_{1}\cdots w_{n-1}}}^{m}d\mu_{-q^{w_{1}\cdots w_{n-1}}}\left(y\right)\frac{t^{m}}{m!}\nonumber \\
= & \sum_{m=0}^{\infty}\left[\prod_{j=1}^{n-1}w_{j}\right]_{q}^{m}\mathcal{E}_{m,q^{w_{1}\cdots w_{n-1}}}\left(w_{n}x+w_{n}\sum_{j=1}^{n-1}\frac{k_{j}}{w_{j}}\right)\frac{t^{m}}{m!}.\nonumber
\end{align}

For $m\ge0$, from (\ref{eq:11}), we have
\begin{align}
 & \int_{\Zp}\left[\left(\prod_{j=1}^{n-1}w_{j}\right)y+\left(\prod_{j=1}^{n}w_{j}\right)x+w_{n}\sum_{j=1}^{n-1}\left(\prod_{\substack{i=1\\
i\neq j
}
}^{n-1}w_{i}\right)k_{j}\right]_{q}^{m}d\mu_{-q^{w_{1}\cdots w_{n-1}}}\left(y\right)\label{eq:12}\\
= & \left[\prod_{j=1}^{n-1}w_{j}\right]_{q}^{m}\mathcal{E}_{m,q^{w_{1}\cdots w_{n-1}}}\left(w_{n}x+w_{n}\sum_{j=1}^{n-1}\frac{k_{j}}{w_{j}}\right),\quad\left(n\in\mathbb{N}\right).\nonumber
\end{align}

Therefore, by Theorem \ref{thm:1} and (\ref{eq:12}), we obtain the
following theorem.
\begin{thm}
\label{thm:2} Let $w_{1},\dots w_{n}\in\mathbb{N}$ be such that $w_{1}\equiv w_{2}\equiv\cdots\equiv w_{n}\equiv1\pmod{2}$.
For $m\ge0$, the following expressions
\begin{align*}
&\frac{\left[\prod_{j=1}^{n-1}w_{\sigma\left(j\right)}\right]_{q}^{m}}{\left[2\right]_{q^{w_{\sigma\left(1\right)}\cdots w_{\sigma\left(n-1\right)}}}}\prod_{l=1}^{n-1}\sum_{k_{l}=0}^{w_{\sigma\left(l\right)}-1}\left(-1\right)^{\sum_{i=1}^{n-1}k_{i}}q^{w_{\sigma\left(n\right)}\sum_{j=1}^{n-1}\left(\prod_{\substack{i=1\\
i\neq j
}
}^{n-1}w_{\sigma\left(i\right)}\right)k_{j}}\\
&\times\mathcal{E}_{m,q^{w_{\sigma\left(1\right)}\cdots w_{\sigma\left(n-1\right)}}}\left(w_{\sigma\left(n\right)}x+w_{\sigma\left(n\right)}\sum_{j=1}^{m-1}\frac{k_{j}}{w_{\sigma\left(j\right)}}\right)
\end{align*}
are the same for any $\sigma\in S_{n}$.
\end{thm}
It is not difficult to show that
\begin{align}
 & \left[y+w_{n}x+w_{n}\sum_{j=0}^{n-1}\frac{k_{j}}{w_{j}}\right]_{q^{w_{1}\cdots w_{n-1}}}\label{eq:13}\\
= & \frac{\left[w_{n}\right]_{q}}{\left[\prod_{j=1}^{n-1}w_{j}\right]_{q}}\left[\sum_{j=1}^{n-1}\left(\prod_{\substack{i=1\\
i\neq j
}
}^{n-1}w_{i}\right)k_{j}\right]_{q^{w_{n}}}+q^{w_{n}\sum_{j=1}^{n-1}\left(\prod_{\substack{i=1\\
i\neq j
}
}^{n-1}w_{i}\right)k_{j}}\left[y+w_{n}x\right]_{q^{w_{1}\cdots w_{n-1}}}.\nonumber
\end{align}

Thus, by (\ref{eq:13}), we get
\begin{align}
 & \int_{\Zp}\left[y+w_{n}x+w_{n}\sum_{j=0}^{n-1}\frac{k_{j}}{w_{j}}\right]_{q^{w_{1}\cdots w_{n-1}}}^{m}d\mu_{q^{-w_{1}\cdots w_{n-1}}}\left(y\right)\label{eq:14}\\
= & \sum_{l=0}^{m}\binom{m}{l}\left(\frac{\left[w_{n}\right]_{q}}{\left[\prod_{j=1}^{n-1}w_{j}\right]_{q}}\right)^{m-l}\left[\sum_{j=1}^{n-1}\left(\prod_{\substack{i=1\\
i\neq j
}
}^{n-1}w_{i}\right)k_{j}\right]_{q^{w_{n}}}^{m-l}q^{lw_{n}\sum_{j=1}^{n-1}\left(\prod_{\substack{i=1\\
i\neq j
}
}^{n-1}w_{i}\right)k_{j}}\nonumber\\
&\times\int_{\Zp}\left[y+w_{n}x\right]_{q^{w_{1}\cdots w_{n-1}}}^{l}d\mu_{-q^{w_{1}\cdots w_{n-1}}}\left(y\right)\nonumber \\
= & \sum_{l=0}^{m}\binom{m}{l}\left(\frac{\left[w_{n}\right]_{q}}{\left[\prod_{j=1}^{n-1}w_{j}\right]_{q}}\right)^{m-l}\left[\sum_{j=1}^{n-1}\left(\prod_{\substack{i=1\\
i\neq j
}
}^{n-1}w_{i}\right)k_{j}\right]_{q^{w_{n}}}^{m-l}\nonumber\\
&\times q^{lw_{n}\sum_{j=1}^{n-1}\left(\prod_{\substack{i=1\\
i\neq j
}
}^{n-1}w_{i}\right)k_{j}}\mathcal{E}_{l,q^{w_{1}\cdots w_{n-1}}}\left(w_{n}x\right).\nonumber
\end{align}

From (\ref{eq:14}), we have
\begin{align}
 &\relphantom{=} \frac{\left[\prod_{j=1}^{n-1}w_{j}\right]_{q}^{m}}{\left[2\right]_{q^{w_{1}\cdots w_{n-1}}}}\prod_{l=1}^{n-1}\sum_{k_{l}=0}^{w_{l}-1}\left(-1\right)^{\sum_{l=1}^{n-1}k_{l}}q^{w_{n}\sum_{j=1}^{n-1}\left(\prod_{\substack{i=1\\
i\neq j
}
}^{n-1}w_{i}\right)k_{j}}\nonumber\\
&\relphantom{=}\times\int_{\Zp}\left[y+w_{n}x+w_{n}\sum_{j=1}^{n-1}\frac{k_{j}}{w_{j}}\right]_{q^{w_{1}\cdots w_{n-1}}}^{n}d\mu_{-q^{w_{1}\cdots w_{n-1}}}\left(y\right)\label{eq:15}\\
 &= \sum_{l=0}^{m}\binom{m}{l}\frac{\left[\prod_{j=1}^{n-1}w_{j}\right]_{q}^{l}}{\left[2\right]_{q^{w_{1}\cdots w_{n-1}}}}\left[w_{n}\right]_{q}^{m-l}\mathcal{E}_{l,q^{w_{1}\cdots w_{n-1}}}\left(w_{n}x\right)\nonumber\\
&\relphantom{=}\times\prod_{s=1}^{n-1}\sum_{k_{s}=0}^{w_{s}-1}\left(-1\right)^{\sum_{j=1}^{n-1}k_{j}}q^{\left(l+1\right)w_{n}\sum_{j=1}^{n-1}\left(\prod_{\substack{i=1\\
i\neq j
}
}^{n-1}w_{i}\right)k_{j}}\left[\sum_{j=1}^{n-1}\left(\prod_{\substack{i=1\\
i\neq j
}
}^{n-1}w_{i}\right)k_{j}\right]_{q^{w_{n}}}^{m-l}\nonumber \\
 &= \frac{1}{\left[2\right]_{q^{w_{1}w_{2}\cdots w_{n-1}}}}\sum_{l=0}^{m}\binom{m}{l}\left[\prod_{j=1}^{n-1}w_{j}\right]_{q}^{l}\left[w_{n}\right]_{q}^{m-l}\mathcal{E}_{l,q^{w_{1}\cdots w_{n-1}}}\left(w_{n}x\right)\nonumber\\
&\relphantom{=}\times\hat{T}_{m,q^{w_{n}}}\left(w_{1},w_{2},\dots,w_{n-1}\mid l\right),\nonumber
\end{align}
where
\begin{align}
&\relphantom{=}\hat{T}_{m,q}\left(w_{1},\dots,w_{n-1}\mid l\right)\label{eq:16}\\
&=\prod_{s=1}^{n-1}\sum_{k_{s}=0}^{w_{s}-1}q^{\left(l+1\right)\sum_{j=1}^{n-1}\left(\prod_{\substack{i=1\\
i\neq j
}
}^{n-1}w_{i}\right)k_{j}}\left[\sum_{j=1}^{n-1}\left(\prod_{\substack{i=1\\
i\neq j
}
}^{n-1}w_{i}\right)k_{j}\right]_{q}^{m-l}\left(-1\right)^{\sum_{j=1}^{n-1}k_{j}}.\nonumber
\end{align}

As this expression is invariant under any permutation in $S_{n}$, we have the following
theorem.
\begin{thm}
\label{thm:3} Let $w_{1},w_{2},\dots,w_{n}\in\mathbb{N}$ be such that
$w_{1}\equiv w_{2}\equiv\cdots\equiv w_{n}\equiv1\pmod{2}$. For $m\ge0$,
the following expressions
\begin{align*}
&\relphantom{=}\frac{1}{\left[2\right]_{q^{w_{\sigma\left(1\right)}w_{\sigma\left(2\right)}\cdots w_{\sigma\left(n-1\right)}}}}\sum_{l=0}^{m}\binom{m}{l}\left[\prod_{j=1}^{n-1}w_{\sigma\left(j\right)}\right]_{q}^{l}\left[w_{\sigma\left(n\right)}\right]_{q}^{m-l}\nonumber\\
&\times\mathcal{E}_{l,q^{w_{\sigma\left(1\right)}\cdots w_{\sigma\left(n-1\right)}}}\left(w_{\sigma\left(n\right)}x\right)\hat{T}_{m,q^{w_{\sigma\left(n\right)}}}\left(w_{\sigma\left(1\right)},w_{\sigma\left(2\right)},\dots,w_{\sigma\left(n-1\right)}\mid l\right)
\end{align*}
are the same for any $\sigma\in S_{n}$.
\end{thm}

Acknowledgements. This paper is supported by grant NO 14-11-00022 of Russian Scientific Fund.

\bibliographystyle{amsplain}
\providecommand{\bysame}{\leavevmode\hbox to3em{\hrulefill}\thinspace}
\providecommand{\MR}{\relax\ifhmode\unskip\space\fi MR }
\providecommand{\MRhref}[2]{%
  \href{http://www.ams.org/mathscinet-getitem?mr=#1}{#2}
}
\providecommand{\href}[2]{#2}

\end{document}